\documentstyle[11pt]{article}
\def\R{\ifmmode{\rm I\mkern-3.1mu
R\mkern1mu}\else{\rm I\kern-.18em  R\hskip1pt\
}\fi\relax}

\def\l{\lambda}

\def\sou{\overline}

\def\f{\rightarrow}

\def\[{\ifmmode{ [\mkern-6.6mu
[\mkern2mu}\else{ [\kern-.28em
[\hskip1pt\ }\fi\relax} 

\def\]{\ifmmode{ ]\mkern-6.6mu
]\mkern2mu}\else{ ]\kern-.28em
]\hskip1pt\ }\fi\relax}

\newtheorem{lemma}{Lemma}[section]

\begin{document}

\begin{center}
\Large\bf
A syntactical proof of the operational equivalence of two
$\l$-terms \\ 
\end{center}

\begin{center}
{\bf Ren\'e DAVID and  Karim NOUR}\\
\end{center}

{\bf Abstract} {\it In this paper we present a purely syntactical proof of the operational
equivalence of $I=\l xx$ and the $\l$-term $J$ that is the $\eta$ -infinite expansion of $I$.}

\section{Introduction}
Two $\l$-terms $M$ and  $N$ are operationnely equivalent ($M \simeq_{oper} N$) iff for all
context $C$ : $C[M]$ is solvable iff $C[N]$ is solvable.

Let $I=\l xx$ and  $J=(Y ~ G)$ where $Y$ is the Turing's fixed point operator and 
$G=\l x \l y \l z (y ~ (x ~ z))$.

$J$ is the $\eta$-infinte expansion of $I$. His B\H{o}hm tree is in fact $\l x \l x_1(x ~ \l
x_2(x_1 ~ \l x_3(x_2 ~ \l x_4(x_3 ...$.

The following Theorem is well known (see [1],[3]).\\

{\bf Theorem} $I \simeq_{oper} J.$ \\

The usual proof is semantic : two $\l$-terms are  operationnely equivalent iff they have the same
interpretation in the modele $D_{\infty}$ . \\

We give below an elementary and a purely syntactical proof of this result. This proof analyses in a
fine way the reductions of $C[I]$ and $C[J]$ by distinguant the "real" $\beta$ -redex of ceux
which come of the $\eta$-expansion.

This proof may be generalize to prove (this result is also well known) the operationnely
equivalence of two $\l$-terms where the B\H{o}hm tree are equal \`a  $\eta$
- infinite expansion  pr\`es. The necessary technical tool is the directed $\l$-calculus (see
[2]).\\

\section{Definitions and notations}

\begin{itemize}
\item $\l \sou{x}~ U$ represents a sequence of abstractions.
\item Let $T,U,U_1,...,U_n$ be $\l$-terms, the application of $T$ to $U$
is denoted by $(T ~ U)$ or $TU$.  In the same way we write $TU_1...U_n$ or $T\sou{U}$
instead of $(...(T ~ U_1)...U_n)$.
\item Let us recall that a $\l$-term $T$ either has a head redex [i.e. $t=\l \sou{x} (\l x U ~ V) ~ \sou{V}$, the head redex being $(\l x U ~ V)$], or is in head
normal form [i.e. $t=\\l \sou{x} x ~ \sou{V}$].
\item The notation $U \f_{t} V$ (resp. $U \f_{t^*} V$) means that $V$ is obtained from $U$ by
one head reduction (resp. some head reductions).  
\item A $\l$-term $T$ is said solvable iff the head reduction of $T$ terminates.
\end{itemize}

The following Lemma is well known.

\begin{lemma} $(U ~V)$ is solvable iff $U$ is solvable (and has $U'$ as head normal form) and  $(U'~
V)$ is solvable.
\end{lemma}

\section {Proof of the Theorem}

The idea of the proof is the following : we prove that, if we assimilate the
reductions where $I$ (resp $J$) are in head position, $C[I]$ and  $C[J]$ reduse, by
head reduction in the same way. For this we add a constante $H$ (which represente
either $I$ or $J$). We define on those terms the $I$ (resp $J$) head reduction,
corresponding to the case where $H = I$ (resp $J$). To prove that the reductions are
equivalent we prove that the terms obtained by "removing" the constante $H$ are equal.
This is the role of the extraction fonction $E$.

\subsection{$\l H$-calculus and the application $E$}

\begin{itemize}
\item We add a new constante $H$ to the $\l$-calculus and we call $\l H$-terms
the terms which we obtain.
\item We define (by induction) on the set of $\l H$-terms the
application $E$ :
\begin{itemize}
\item[] $E(x)=x$ ; $E(H)=H$ ; $E(\l x U)=\l x E(U)$ ; \\
$E(U V)= E(U) E(V)$ if $U \not =H U_1 U_2 ... U_n$ ;
\item[] $E(H U_1 U_2 ... U_n)= E(U_1 U_2 ... U_n)$ .
\end{itemize}
\item A $\l H$-term is in head normal form if it is of the forme
: $\l \sou{x}~ H$ or $\l
\sou{x}~ x \sou{V}$. \end{itemize}

\begin{lemma} If $T$ is a $\l H$-term, then  $E(T)$ is of the forme $\l
\sou{x}~ H$ or $\l \sou{x}~ x\sou{V}$ or $\l \sou{x}~ (\l x U ~ V) \sou{V}$.
\end{lemma}

{\bf Proof} By induction on $T$. $\Box$

\begin{lemma} If $T$ is a $\l H$-term, then  $E(E(T))=E(T)$.
\end{lemma}

{\bf Proof} By induction on $T$. $\Box$

\begin{lemma} Let $T,\sou{U}$ be $\l H$-terms.
$E(T\sou{U})=E(E(T)\sou{E(U)})$.
\end{lemma}

{\bf Proof} By induction on $T$. We distinguish the cases: $T \not = H\sou{V}$
and $T = H\sou{V}$. $\Box$

\begin{lemma} Let $U,V$ be $\l H$-terms and $x$ a
variable, $E(U[V/x])=E(E(U)[E(V)/x])$ . \end{lemma}

{\bf Proof}  By induction on $U$. The only interesting case is $U=x\sou{U}$. By
Lemma 3.3, $E(U[V/x])=E(E(V)\sou{E(U[V/x])})$. Therefore, by induction hypothesis and 
Lemma 3.3,\\ $E(U[V/x])=E(E(V)\sou{E(E(U)[E(V)/x])})=E(E(U[E(V)/x])$. $\Box$

\begin{lemma} Let $U_1,U_2,V_1,V_2$ be $\l H$-terms such that
$E(U_1)=E(U_2)$ and
$E(V_1)=E(V_2)$. $E(U_1[V_1/x])=E(U_2[V_2/x])$.
\end{lemma}

{\bf Proof} By Lemma 3.4. $\Box$

\begin{lemma} Let $U_1,U_2,V_1,V_2$ be $\l H$-terms.
If $U_1 \f_t V_1$, $U_2 \f_t V_2$, and  $E(U_1)=E(U_2)$, then  $E(V_1)=E(V_2)$.
\end{lemma}

{\bf Proof} By Lemmas 3.3 and  3.5. $\Box$

\subsection{The $I$-reduction}

\begin{itemize}
\item We define on the  $\l H$-terms  a new  head reduction :
\begin{itemize}
\item[]  $HU_1...U_n \f_{I} U_1U_2...U_n$
\end{itemize}
\item We denote by $\f_{I^*}$ the reflexive and transitive closure of $\f_{I}$.
\item  A $\l H$-term $U$ is $I$-$t$-solvable iff a finite sequence of $I$-reductions and 
$t$-reductions of $U$ gives a head normal form.

\end{itemize}

\begin{lemma}  Let $U,V$ be $\l H$-terms. If $U \f_{I^*} V$, then $E(U)=E(V)$.
\end{lemma}

{\bf Proof} By induction on the reduction of $U$. $\Box$

\begin{lemma}  Each $I$-reduction is finite.
\end{lemma}

{\bf Proof} The $I$-reduction decreases the complexity of a $\l H$-term.
$\Box$

\begin {lemma} Let $U$ be $\l H$-term. U is $I$-$t$-solvable iff $U[I/H]$ is solvable.
\end{lemma}
{\bf Proof} Immediate. $\Box$

\subsection{The $J$-reduction}

\begin{itemize}
\item We define on the $\l H$-terms a new head reduction :
\begin{itemize}
\item[] $HU_1...U_n \f_{J} U_1(H ~ U_2)U_3...U_n$
\end{itemize}
\item We denote by $\f_{J^*}$ the reflexive and transitive closure of $\f_{J}$.
\item  A $\l H$-term $U$ is $J$-$t$-solvable iff a finite sequence of $J$-reductions and 
$t$-reductions of $U$ gives a head normal form.

\end{itemize}

\begin{lemma}  Let $U,V$ be $\l H$-terms. If $U \f_{J^*}V$, then $E(U)=E(V)$.
\end{lemma}

{\bf Proof} It is enough to do the proof for one step of $J$-reduction. The only interesting case is
$U=(H)U_1U_2\sou{U}$. In this case $U \f_J U_1(H~U_2)\sou{U}$, and, by induction hypothesis,
$E((U_1(H~U_2)\sou{U})=E(V)$, therefore -by  Lemma 3.3- $E(U)=E(V)$. $\Box$

\begin{lemma}  Let $U,V$ be $\l H$-terms. If $U \f_{J^*} V$, then,
for each sequence $\sou{W}=W_1...W_n$, there is a sequence $\sou{W'}=W'_1...W'_n$ such that
$U \sou{W} \f_{J^*} V \sou{W'}$ and for, all $1 \leq k \leq n$, $W'_k \f_{J^*} W_k$.
\end{lemma}

{\bf Proof} By induction on the reduction of $U$. It enough to do the proof for one
step of $J$-reduction. The only interesting case is $U=HU'$ and $\sou{W}=W_1\sou{W'}$.
In this case $V=U'$, $UW_1\sou{W'} \f_J V(H ~ W_1)\sou{W'}$ and  $H W_1 \f_J W_1$. $\Box$

\begin{lemma}  Each $J$-reduction is finite.
\end{lemma}

{\bf Proof} By induction on $U$. The only interesting case is $U=HV_1...V_n$ $(n \geq 2)$. We prove,
by recurrence on $n$, that if the reductions of $V_1, ..., V_n$ are finite, then so is for $U =
H V_1 ... V_n$. $U \f_{J} V_1 (H~ V_2)~ V_3 ...V_n$ and  $ V_1 \f_{J^*} V'_1$. By Lemma 3.11, $U
\f_{J} V'_1 W_2 W_3...W_n$ where $W_2 \f_{J} H~V_2 \f_{J} V_2 $ and  $W_i \f_{J} V_i$, therefore the
reductions of  $W_i$ are finite.

- If $E(V_1) \not = H$.  $V'_1$ begin soit by $\l$, soit by a $\beta$-redex, soit
by a variable. Therefore, by Lemma 3.11, the $J$-reduction of $U$ is finite.\\
- If $E(V_1) = H$. By Lemma 3.11, $U \f_{J^*} HW_2 ...W_n$ and the recurrence hypothesis allows
to conclude. $\Box$

\begin{lemma}  Let $U$ be a $\l H$-term.  $U$ is $J$-$t$-solvable iff
$U[J/H]$) is solvable.  \end{lemma}

{\bf Proof} The only difficulty is to prove that : if $U$ is $J$-$t$-solvable, then $U[J/H]$ is
solvable.\\ We prove that by induction on the reduction of $U$. The only interesting case is $U=\l
\sou{x}~HV$. In this case, $U \f_J \l \sou{x}~ V$ and  $U[J/H]\f_t \l \sou{x}~ \l y
V[J/H]~(J~y)$. By induction hypothesis $V[J/H]$ is solvable, and,
by Lemma 2.1, we may begin to reduse $V[J/H]$ in $\l \sou{x}~ \l y
V[J/H]~(J~y)$. If the head normal form of $V[J/H]$ is not of the forme $\l x
\l\sou{z}~x \sou{W}$, the  result is true. If not the head reduction of $U[J/H]$
gives  $\l \sou{x}~\l\sou{z}~(J~y)\sou{W}$ which is solvable. $\Box$

\subsection{The proof of the Theorem}

$U \f_{(I^*, k)} V$ (resp. $U \f_{(J^*, k)} V$) means that $V$ is obtained from $U$
by $I$-reductions (resp. $J$-reductions) and $k$ $t$-reductions.

\begin{lemma} Let $U_1,U_2,V_1,V_2$ be $\l H$-terms.
If $U_1 \f_{(I^*, k)} V_1$, $U_2 \f_{(J^*, k)} V_2$, and  $E(U_1)=E(U_2)$, then 
$E(V_1)=E(V_2)$. \end{lemma}

{\bf Proof} Consequence of Lemmas 3.6, 3.7 and  3.10. $\Box$

\begin{lemma} Let $U$ be a $\l H$-term. $U$ is $I$-$t$-solvable iff $U$ is
$J$-$t$-solvable.
\end{lemma}

{\bf Proof} Consequence of Lemmas 3.8, 3.12 and 3.14. $\Box$ \\

{\bf Proof of the Theorem } Consequence of Lemmas 3.9, 3.13 and 3.15. $\Box$

LAMA - Equipe de Logique - Universit\'e de Chamb\'ery - 73376 Le Bourget du
Lac \\
e-mail david,nour@univ-savoie.fr

\end{document}